\documentstyle[12pt]{article}

\newcommand{\be}{\begin{equation}}
\newcommand{\ee}{\end{equation}}
\newcommand{\ba}{\begin{eqnarray}}
\newcommand{\ea}{\end{eqnarray}}
\newcommand{\baa}{\begin{eqnarray*}}
\newcommand{\eaa}{\end{eqnarray*}}
\newcommand{\bb}{\begin{thebibliography}}
\newcommand{\eb}{\end{thebibliography}}

\newcommand{\bi}[1]{\bibitem{#1}}
\newcommand{\lab}[1]{\label{#1}}
\newcommand{\re}[1]{(\ref{#1})}

\newcounter{my}
\newcommand{\he}%
   {\stepcounter{equation}\setcounter{my}%
   {\value{equation}}\setcounter{equation}0%
   }%
\newcommand{\she}%
   {\setcounter{equation}{\value{my}}%
    }%

\begin{document}


\vspace*{10mm}

\begin{center}

{\Large \bf Zeros and orthogonality of the Askey-Wilson \\
          polynomials for $q$ a root of unity}%
\footnote{{\em Date:} April 1996; published in {\em Duke Math. J.} {\bf 89} (1997), 283--305.}

\vspace{5mm}

{\large \bf Vyacheslav Spiridonov and Alexei Zhedanov}

\end{center}

\begin{abstract}

We study some properties of the Askey-Wilson polynomials (AWP) when $q$ is a
primitive $N$th root of unity. For general four-parameter AWP,
zeros of the $N$th polynomial and the orthogonality measure
are found explicitly.
Special subclasses of the AWP, e.g., the continuous
$q$-Jacobi and big $q$-Jacobi polynomials, are considered
in detail. A set of discrete weight functions positive on
a real interval is described.
Some new trigonometric identities related to the AWP are obtained.
Normalization conditions of some polynomials are expressed in terms
of the Gauss sums.

\end{abstract}

\bigskip


\medskip


\bigskip


\section{Introduction}
The theory of orthogonal polynomials of one variable is based upon the
three-term recurrence relation \cite{Chi}
\be
P_{n+1}(x) +b_nP_n(x)+u_nP_{n-1}(x)=xP_n(x), \qquad n=1, 2, \dots
\lab{rec1} \ee
with the initial conditions
\be
P_0(x)=1, \quad P_1(x)=x-b_0.
\lab{init1} \ee
There is a representation theorem concerning orthogonality of the
polynomials $P_n(x)$ with respect to
some complex measure $\mu(x)$ in case $u_n$ and $b_n$
are arbitrary complex coefficients; see, e.g., \cite{AW,Chi},
\be
\int_C P_n(x)P_m(x) d\mu(x) = h_n \delta_{nm},
\lab{genort} \ee
where $h_n=u_1u_2\cdots u_n$
and $C$ is some contour in the complex plane. When $b_n$
are real and $u_n>0$, the measure can be chosen real-valued
(then $C$ is the real line) and
nondecreasing. For such cases it can be shown that all zeros $x_s,
s=0,1,\dots,N-1,$ of any polynomial $P_N(x)$ are real and simple \cite{Chi}.
Moreover, these zeros can be used for construction of a
discrete orthogonality relation for polynomials of the degree lower than $N$:
\be
\sum_{s=0}^{N-1} P_n(x_s)P_m(x_s)w_s = h_n \delta_{nm},
\lab{finort} \ee
where $n, m = 0, 1, \dots, N-1$ and the weight function is
\be
w_s =\frac{h_{N-1}}{P_{N-1}(x_s)P_N^{\prime}(x_s)}.
\lab{ws} \ee
The property \re{finort} was proven in \cite{At} under the assumption that
$b_n$ are real and $u_n$ are positive. However, the proof remains valid
even in general cases of complex $b_n, u_n$. The only requirement to be
added is that all zeros $x_s$ must be simple which is not guaranteed if
the recurrence coefficients are complex.
For example, for $b_n=0$ one finds from the recurrence
relation \re{rec1}, $P_3(x)=x(x^2-u_1-u_2)$. If $u_2=-u_1$ then all
zeros of the polynomial $P_3(x)$ are located at $x=0$. It is thus always
necessary to check
simplicity of the zeros $x_s$ in order for \re{finort} to be applicable.

In this paper we present a study of the
finite orthogonality property \re{finort}
of the Askey-Wilson polynomials (AWP) for the case when $q$ is a root of unity.
Recall that the recurrence coefficients for AWP  have the form  \cite{AW,GR}
\be
u_n=\xi_n \eta_{n-1}, \quad b_n=-\xi_n -\eta_n +a +a^{-1},
\lab{recaw} \ee
where
\be
\xi_n= \frac{a (1-q^n)(1-bcq^{n-1})(1-bdq^{n-1})(1-cdq^{n-1})}
{(1-gq^{2n-2})(1-gq^{2n-1})},
\lab{xi} \ee
\be
\eta_n= \frac{a^{-1}(1-abq^n)(1-acq^n)(1-adq^n)(1-gq^{n-1})}
{(1-gq^{2n})(1-gq^{2n-1})},
\lab{eta} \ee
and $g=abcd$. The parameters $a, b, c, d$ are arbitrary complex numbers.
The explicit expression of AWP in terms of the basic hypergeometric
functions is \cite{AW,GR}
\be
P_n(x)= D{_n}\;{_4}\varphi_3 \left({q^{-n},gq^{n-1},at,a/t
\atop ab, ac, ad} ;q, q \right),
\lab{aw} \ee
where
\be
D_n= \frac{a^{-n}(ab;q)_n(ac;q)_n(ad;q)_n}{(gq^{n-1};q)_n}
\lab{C_n} \ee
and $x=t+t^{-1}$. Here $(a;q)_n$ is the $q$-shifted factorial defined
as $(a;q)_0=1$ and $(a;q)_n=(1-a)(1-aq)\cdots(1-aq^{n-1})$ for integer $n>0$.
The general $q$-hypergeometric series is defined as follows \cite{GR}:
\be
{_r}\varphi_s\left({a_1, a_2, \dots, a_r \atop b_1, \dots, b_s}; q, z \right)=
\sum_{n=0}^\infty{(a_1, a_2, \dots, a_r; q)_n \over
(q, b_1, \dots, b_s; q)_n} [(-1)^nq^{n(n-1)/2}]^{1+s-r} z^n,
\lab{qser}\ee
where the compact notation $(a_1, a_2, \dots, a_r; q)_n =
(a_1; q)_n\cdots (a_r; q)_n$ is used.

For real $q$ and appropriate restrictions upon the parameters $a,b,c,d$
the AWP are orthogonal on a finite interval of the
real axis with a positive continuous weight function \cite{AW}.
In this case, AWP satisfy a difference equation defined upon the
hyperbolic grid \cite{AW,NSU}.
A possibility of constructing orthogonal polynomials for trigonometric
grids has been mentioned briefly in \cite{NSU}. In this paper
we concentrate on the case when $q$ is a primitive $N$-th root of unity:
\be
q=\exp(2{\pi}iM/N),
\lab{q} \ee
where $M$ and $N$ are mutually co-prime positive integers, $M<N$.
We consider only finite sets of orthogonal polynomials with discrete
measures which are
different from the so-called sieved polynomials (see, e.g., \cite{AAA}).

\section{Zeros of the $N$th Askey-Wilson polynomial}
\setcounter{equation} 0
In this section we analyze location of zeros for generic $N$th order
AWP depending on four parameters for $q$ a primitive $N$th root of unity.

Consider the explicit form of $P_N(x)$ \re{aw} for $q$ fixed by
\re{q}. Obviously we cannot use this expression directly, because the
factor $(q^{-N};q)_N/(q;q)_N$ in the $N$th term of the series is
indeterminate. Nevertheless, we can apply a
limiting procedure taking $q=\exp(\epsilon + 2\pi iM/N)$ and going to the
limit $\epsilon \to 0$. Then it is easy to see that
\be
\lim_{\epsilon \to 0}\; {(q^{-N};q)_N \over (q;q)_N} =
\lim_{\epsilon\to 0}\; (-1)^N q^{-N(N+1)/2}=-1.
\lab{lim} \ee
Using \re{lim} one can find that
\be
P_N(x)=D_N \left(1-\frac{(g;q)_N(at;q)_N (at^{-1};q)_N }
{(ab;q)_N (ac;q)_N (ad;q)_N}\right).
\lab{PN} \ee
Since $u_N=0$ in \re{rec1}. we have a finite set of polynomials
and it is not necessary to consider ambiguities in the higher order terms
in the series \re{aw}.

The $q$-shifted factorial in \re{PN} can be simplified because for any
complex parameter $a$ one has the identity (which is valid only for $q$ a root
of unity)
\be
(a;q)_N=1-a^N.
\lab{Npoch} \ee
Indeed, $(a;q)_N=(1-a)(1-aq)\cdots(1-aq^{N-1})$ is a $N$-th order polynomial
in argument $a$. The zeros
of this polynomial are $1,q,q^2, \dots, q^{N-1},$ whence \re{Npoch} follows.

We can rewrite \re{PN} in the form
\be
P_N(x)=t^N + t^{-N} - E_N,
\lab{PN2} \ee
where $E_N$ is a rational function of parameters,
\be
E_N=
{a^N+b^N+c^N+d^N -(abc)^N-(bcd)^N-(abd)^N-(acd)^N \over 1-(abcd)^N}.
\lab{EN} \ee
We call $E_N$ a {\it combinatorial invariant} because its symmetry
with respect to permutation of the parameters $a, b, c, d$  plays an
important role in the following.

From \re{PN2} it is not difficult to find zeros of $P_N(x)$:
\be
x_s= rq^s+r^{-1}q^{-s}, \qquad
s=0,1,\dots,N-1, \lab{roots} \ee
where for simplicity we assume that $r$ is the root with minimal argument
of the equation
$$
r^N= E_N/2 +\sqrt{E_N^2/4 -1}, \quad \mbox{\rm or} \quad
r^{-N} = E_N/2-\sqrt{E_N^2/4-1}.
$$
Let us point out a curious identity associated with the combinatorial
invariant $E_N$:
\be
\frac{(1-a^Nr^N)(1-b^Nr^N)(1-c^Nr^N)(1-d^Nr^N)}
{(a^N-r^N)(b^N-r^N)(c^N-r^N)(d^N-r^N)}=1.
\lab{identity}\ee

Now we are ready to formulate the orthogonality relation for AWP. First
we need the condition of simplicity of the zeros $x_s$, which is
\begin{equation}
E_N \ne \pm 2.
\label{encond}\end{equation}
Second, in order for the polynomials $P_n(x), 0\leq n \leq N,$ to be well
defined we impose the constraints
\be
g, ab, ac, ad, bc, bd, cd  \ne q^k, \quad k=0, 1, \dots, N-1,
\lab{cond1} \ee
where $q$ is fixed by \re{q}.
These restrictions guarantee that $u_n$ does not vanish between two zeros
at $n=0$ and $n=N$ or that $u_n, b_n$ do not blow up (cf. \cite{Sz1}).
Under the conditions \re{cond1} one can calculate
\be
P_N^{\prime}(x_s)= N\frac{r^N-r^{-N}}{rq^s-r^{-1}q^{-s}}.
\lab{P'}\ee
Using the same $\epsilon\to 0$ limiting procedure as in \re{lim}
for (unambiguous) determination of the polynomial $P_{N-1}(x_s)$ we get
the following

{\bf Theorem 1.}
{\it The weight function of the
finite-dimensional set of Askey-Wilson polynomials for $q$ a primitive root
of unity \re{q} and the constraints upon the parameters
\re{encond}, \re{cond1} has the form
\be
w_s=  \frac{h_{N-1}(rq^s-r^{-1}q^{-s})}{ND_{N-1}(r^N-r^{-N})F(s)},
\lab{w1} \ee
where $s=0, 1, \dots, N-1, $ and
\be
F(s)=\sum_{n=0}^{N-1} \frac{(gq^{-2}, arq^s, ar^{-1}q^{-s}; q)_n}
{(ab, ac, ad; q)_n}q^n.
\lab{F(s)} \ee
}
\medskip

This form of the weight function is not convenient for
the analysis of its positivity since a truncated ${_4} \varphi_3$ series
\re{F(s)} depending on the discrete argument $s$ is involved
into the definition. In the next section we shall derive a more explicit
form of $w_s$.

\section{Difference equation and absence of self-duality}
\setcounter{equation} 0

In this section we derive a difference equation for generic AWP for $q$ a
root of unity and show that a complete self-duality property
(a permutational symmetry between the parametrization of
polynomial's argument and its degree) is absent.

In order to derive this difference equation we use contiguous relations for
the terminating basic hypergeometric series ${_4}\varphi_3$.
Recall that if one denotes
$$
\phi(\xi, \eta)= {_4} \varphi_3 \left({\xi,\eta,\gamma,\delta
\atop e, f, h} ;q, q \right),
$$
where one of the parameters $\xi,\eta,\gamma,\delta$ is equal to $q^{-M}$,
$M$ a positive integer, and $efh=\xi\eta\gamma\delta q$,
then one has the following relation \cite{AW0,GR}:
\be
B_1\phi(\xi q^{-1}, \eta q) + B_3 \phi(\xi, \eta) + B_2 \phi(\xi q , \eta
q^{-1})=0, \lab{contig} \ee
where
\ba
B_1&=&\eta(1-\eta)(\xi q-\eta)(\xi-e)(\xi-f)(\xi-h),\\
B_2&=&\xi(1-\xi)(\xi-\eta q)(\eta -e)(\eta-f)(\eta -h),\\
B_3&=&\xi \eta (\xi -\eta q)(\xi q - \eta)(\xi -\eta)(1-\gamma)(1-\delta)
- B_1 - B_2.  \ea
In our case one should choose the parametrization
\be \xi=arq^s, \; \eta =
ar^{-1}q^{-s},\; \gamma = q^{-n}, \delta = gq^{n-1}, \; e=ab, \; f=ac, \;
h=ad, \lab{conpar} \ee
which gives $P_n(x_s)=D_n\phi(arq^s, ar^{-1}q^{-s}).$
Substituting \re{conpar} into \re{contig}
we arrive at the difference equation for AWP
with the argument $x_s $ given by \re{roots}:
\be A_sP_n(x_{s-1})
+C_sP_n(x_{s+1}) - (A_s +C_s)P_n(x_s)= \lambda_n P_n(x_s),
\lab{diff} \ee
where
\be
\lambda_n = (q^{-n}-1)(1-gq^{n-1})
\lab{lamb} \ee
and
\ba
A_s&=&\frac{gq^{-1}(1-ra^{-1}q^s)(1-rb^{-1}q^s)(1-rc^{-1}q^s)
(1-rd^{-1}q^s)}{(1-r^2q^{2s-1})(1-r^2q^{2s})},\\
C_s&=&\frac{(1-arq^s)(1-brq^s)(1-crq^s)(1-drq^s)}
{(1-r^2q^{2s})(1-r^2q^{2s+1})}.
\lab{AC} \ea
Let us note that after the (complex) shift of the variable $s$ such
that $ra^{-1}q^s\to q^s$ and a change of notation for parameters,
the equation \re{diff} coincides with the three-term recurrence relation for
AWP. This means that the difference equation \re{diff} may be considered as
a recurrence relation for a special class of
associated AWP. The general associated AWP are determined by \re{rec1},
\re{init1} and the coefficients $u_n, b_n$ obtained from \re{recaw}
after the replacement of $q^n$ by $\mu q^n$, where $\mu$ is an
arbitrary complex parameter \cite{IR}.

But the recurrence coefficients $A_s$ and $C_s$ do not have zeros for any
$s=0, 1, \dots,$ so that one can consider orthogonal polynomials
of unlimited degree. The orthogonality measure of these polynomials
will be continuous, i.e. the parameter $\lambda_n$ on the right
hand side of \re{diff} will not be quantized. Moreover, there exist
regions of the parameters $a, b, c, d$ such that the coefficients
$u_n$ are positive and $b_n$ are real. For instance, such a situation takes
place if we set
\be
a=q^\alpha, \quad d=q^\delta, \quad b=-q^\beta, \quad c=-q^\gamma,
\quad q=e^{2\pi i/N},
\lab{positpar}\ee
where $\alpha, \beta, \gamma, \delta$ are real nonzero parameters,
choose $N$ even, and impose the restrictions
\be
|\alpha|, |\beta|, |\gamma|, |\delta|\ll 1,
\quad \alpha<0,\quad \alpha<\beta<-\alpha,\quad \gamma>-\alpha,
\quad \delta>-\alpha.
\lab{restric}\ee
Under these conditions the difference equation \re{diff} becomes
Hermitian. From this fact follows the existence of $q$-orthogonal
polynomials (a special subclass of the associated AWP) for $q$
a root of unity with the continuous positive measure which are
different from the sieved polynomials. Consideration of these
polynomials lies beyond the scope of present work.

Before calculating the weight function, let us discuss a
self-duality of the AWP \cite{BI,Leo}. For simplicity let us limit
consideration to finite sets of orthogonal polynomials with positive measures
$P_n(s)$, where $n$ is the degree of polynomials and $s$ is the index
numbering discrete values of the argument $x_s$,
$n, s= 0, 1, \dots, N-1>0.$
Due to the Favard theorem, any set of such polynomials is
uniquely determined by the three-term recurrence relation.
As usual, one has some discrete orthogonality relations \re{finort} and
dual conditions:
\be
\sum_{n=0}^{N-1} \frac{P_n(s)P_n(s')}{h_n}=w_s^{-1}\delta_{ss'}
\lab{dualorth}\ee
(see, e.g., \cite{At}).
In the dual picture,$P_n(s)$ are considered as eigenvectors with
eigenvalues $\lambda_n$ of some Hermitian matrix acting upon the
(dual) index $s$. The latter matrix is not necessarily of the Jacobi form;
in general it is difficult to find explicitly. Moreover,
it is not unique (there is a dependence on functional parametrization
of $x_s$ on $s$). When this matrix is tridiagonal, one may
refer to the corresponding polynomials as partially self-dual,
because then both indices $n$ and $s$ can be considered as degrees
of orthogonal polynomials. Polynomials
$P_n(s)$ are called completely self-dual if there exists a symmetry such
that $P_n(s)=P_s^*(n)$ (this relation requires a special normalization of
polynomials, see \cite{BI,Leo}). Here the involution sign $'{^*}'$ means a
permutation of parameters of polynomials.
In the simplest case there is no transformation of parameters at all and the
symmetry condition is $P_n(s)=P_s(n)$, which
means that $P_n(s)$, being considered as a $N\times N$ matrix, equals
 its transpose. For completely self-dual systems,
the three-term recurrence relation in $n$ and the difference equation in
$s$ coincide up to a change of parameters. In the $P_n(s)=P_s(n)$
case, these equations just totally coincide. For examples
of such a situation we refer to the discrete series representations of
the quantum algebra $so_q(3)$ \cite{Sz2} describing
$q$-ultraspherical polynomials at $q^N=1$.

In \cite{BI,Leo} completely self-dual orthogonal polynomials were considered
and it was proven that the only
polynomials possessing this property are the $q$-Racah polynomials.
The recurrence coefficients $\xi_n, \: \eta_n$ in \re{recaw} for
these polynomials \cite{AW0} satisfy the conditions
$\xi_0=\eta_{N+1}=0$.
For the case $q^N=1$ with generic values of parameters
only the first factor $\xi_n$ creates zeros $\xi_0=\xi_{N}=0$.
I.e., our $P_n(x_s)$ do not belong to the
$q$-Racah polynomials and hence they are not completely
self-dual. The latter can also be seen directly
from the explicit expression \re{aw}:
for $t=rq^s$ the permutation of $n$ and $s$ cannot be compensated
by an appropriate normalization of polynomials and permutation of
parameters. In fact, the values of variable $s$ can be
extended from a finite set $s=0,1,\dots, N-1$ to all integers without any
problem. As a result, one gets the periodicity property
$P_n(x_{s+N})=P_n(x_s)$ and the
same periodicity is satisfied by the coefficients $A_s,
C_s$ in the difference equation \re{diff}. Evidently, these
$A_s, C_s$ never vanish. $A_sC_{s-1} \neq 0$, so that there is
not automatic reduction of the difference equation to the
finite-dimensional matrix form in contrast to the $q$-Racah polynomials case.

As seen from \re{aw}, the polynomials $P_n(x_s)$ become completely self-dual
if the condition $r=aq^j$, $j$ integer, is fulfilled. But this equality
holds only if one of the products $ab,\, ac,\, ad$ is equal to $q^k$,
$k$ integer, as it follows from \re{EN}. However, such constraints
correspond to the $q$-Racah polynomials and they are forbidden
in our case \re{cond1}. To summarize the above
considerations: duals of the AWP for $q^N=1$ are given by the special
types of the associated AWP for $q^N=1$, i.e., there is no complete
self-duality.

Let us return now to the problem of calculation of the weight function $w_s$.
Existence of the difference equation \re{diff} allows us to find a simpler
form of $w_s$ than \re{w1}.
Indeed, multiplying equation \re{diff} by a function $\sigma(s)$ satisfying
the only requirement of periodicity $\sigma(s+N)=\sigma(s)$,
and combining it with a similar equation for the polynomials $P_m(x_s)$,
one can get a bilinear relation:
\ba \nonumber
A_s\sigma(s)[P_n(x_{s-1})P_m(x_s)-P_n(x_s)P_m(x_{s-1})] \\ \nonumber
+ C_s\sigma(s)[P_n(x_{s+1})P_m(x_s)-P_n(x_s)P_m(x_{s+1})] \\
= (\lambda_n-\lambda_m) \sigma(s) P_n(x_s)P_m(x_s).
\lab{diff1} \ea
Choose $\sigma(s)$ in such a way that
\be
A_{s+1}\sigma(s+1)=C_s\sigma(s). \lab{sigma}
\ee Then
the equation \re{diff1} can be rewritten in the form
\be
R_{nm}(s+1)-R_{nm}(s)= (\lambda_n - \lambda_m) \sigma(s) P_n(x_s) P_m(x_s)
\lab{diff2}, \ee
where
\be
R_{nm}(s)= A_s\sigma(s)[P_n(x_s)P_m(x_{s-1})-P_n(x_{s-1})P_m(x_s)].
\lab{Rnm} \ee
Summing \re{diff2} from $s=0$ to $s=N-1$ and using the obvious periodicity
property $R_{nm}(s+N)=R_{nm}(s)$, one gets the orthogonality relation:
\be
(\lambda_n - \lambda_m) \sum_{s=0}^{N-1} {\sigma(s) P_n(x_s) P_m(x_s)} = 0,
\quad n \ne m. \lab{orth2} \ee
The weight function $\sigma(s)$ is determined from the relation \re{sigma}
uniquely, up to a multiplicative factor. Hence one should have
$w_s=\kappa \sigma(s)$ where $\kappa$ is some constant and $w(s)$ is the
weight function calculated in \re{w1}. From \re{sigma} we find:
\be w_s= w_0 (q/g)^s
{(1-r^2q^{2s})(ar, br, cr, dr; q)_s\over (1-r^2)(qr/a, qr/b, qr/c, qr/d; q)_s},
\lab{wdiff} \ee
where $w_0$ is a normalization constant. Note that the consistency condition
$\sigma(N)=\sigma(0)$ imposed upon $q$-difference equation \re{sigma} for
$q^N=1$ (see \cite{Car} for a general discussion of such conditions)
is satisfied in a quite nontrivial manner due to \re{identity}.

We have found thus the explicit dependence of the weight function $w_s$
on the discrete argument $s$. The remaining value $w_0$ can be evaluated
from \re{w1}
\be
w_0=a^{N-1}\frac{(1-(bc)^N)(1-(cd)^N)(1-(bd)^N)(1-g/q)
(1-g/q^2)(r-r^{-1})}{F(0)(1-g^N)^2 (1-bc/q)(1-cd/q)(1-bd/q)(r^N-r^{-N})},
\lab{w(0)} \ee
where $F(0)$ is fixed by \re{F(s)} at $s=0$. The normalization condition
\be
\sum_{s=0}^{N-1} w_s =1,
\lab{gennorm} \ee
where $w_s$ is given by \re{wdiff}, leads to some non-trivial trigonometric
identities which we summarize into the following theorem.

{\bf Theorem 2.} {\it Under the constraints \re{encond}, \re{cond1} upon the
complex parameters $a, b, c, d$ and for $q=e^{2\pi iM/N}$, $(N, M)=1$,
the following trigonometric identity takes place
\be
\sum_{s=0}^{N-1}\left({q\over g}\right)^s\frac{(1-r^2q^{2s})
(ar, br, cr, dr; q)_s}{(qr/a, qr/b, qr/c, qr/d; q)_s}
\lab{awpid} \ee
$$
=\frac{r^{1-N}(1-r^{2N})(1-g^N)^2 f(a, b, c, d)}
{(ab, ac, ad, bc, cd, bd; q)_{N-1} (1-g/q)(1-g/q^2)},
$$
where
\be
f(a, b, c, d)\equiv \frac{(ab, ac, ad; q)_{N-1}}{a^{N-1}}
\sum_{n=0}^{N-1} \frac{(gq^{-2}, ar, ar^{-1}; q)_n}
{(ab, ac, ad; q)_n}q^n.
\lab{fabcd}\ee
}

\medskip

{\it Comments.} Probably the truncated $_4\varphi_3$ series \re{fabcd}
has a simpler representation through some finite products due to the
specific properties of the roots of unity and algebraic relation
between $r$ and $a, b, c, d$,
but we were not able to find it. Since the AWP are symmetric
with respect to permutation of parameters, the same
symmetry should be valid for $f(a,b, c, d)$, but it has not been
proved yet to be explicit. Since the left hand side of \re{awpid}
is finite for $a\to 0$, the function \re{fabcd} should be finite
in this limit as well. The above identity looks as a root-of-unity
analog of a special case of the Watson transformation formula for
$_8\varphi_7$ series \cite{GR}, but the complete coincidence is not
yet established.

The relation \re{awpid} exists only for $q$ a root of unity,
i.e., there are no its analogs for other values of $q$ contrary
to the $q$-Racah polynomials' weight function normalization condition
holding for $|q|<1, |q|>1$, or $|q|=1$. As shown below (see also the
Appendix), at the bottom of trigonometric identities obtained for $q^N=1$
like \re{awpid}, one finds Gauss sums or similar summation
formulas. The latter play an important role in the number theory.
Therefore, the arithmetic aspects of \re{awpid} require
a detailed investigation. Note in passing that the most general identity
of such type should be tied to the biorthogonal rational functions
associated with the special $_{10}\varphi_9$ series \cite{mas}.

\section{Two parameter symmetric Askey-Wilson polynomials}
\setcounter{equation} 0

In this section we consider a special subclass of the symmetric
AWP, and determine some regions of parameters for which the measure
is positive.

Consider AWP with $a=-c=q^{\alpha}, \: d=-b=q^{\beta}$.
Then in the recurrence relation for AWP only nondiagonal terms survive,
i.e., the polynomials are symmetric:
\be
P_{n+1}(x) +u_n P_{n-1}(x) = x P_n(x),
\lab{sym} \ee
where $n =0,1,\dots,N-1$,
\be
u_n=\frac{4\sin\omega n \, \sin \omega(n+ 2\alpha + 2\beta -2) \, \cos
\omega(n+2\alpha-1) \, \cos \omega(n +2\beta-1)}{\sin 2\omega(n+\alpha
+\beta -3/2) \, \sin 2\omega(n +\alpha + \beta -1/2)}
\lab{usym} \ee
and $\omega =\pi M/N$.

One can show that if $M=1$ and $N$ is odd, the recurrence coefficients
$u_n$ are positive, $u_n>0,$ for $n= 1, 2, \dots, N-1,$ under the following
constraints upon the parameters $\alpha$ and $\beta$:
\be
-1/4 < \alpha < 1/4, \qquad 1/4 < \beta < 3/4,
\lab{I} \ee
or
\be
1/4 < \alpha < 3/4, \qquad 3/4 < \beta < 5/4.
\lab{II} \ee
These restrictions, or their partners obtained after permuting the
parameters $\alpha$ and $\beta$, provide only particular examples of
positive $u_n$, the problem of finding all such cases (including the cases
of $M \ne 1$ and even $N$) is complicated and is not considered here.

The difference equation for the taken symmetric AWP has the form \re{diff},
where
\ba
A_s &=& \frac{\sin 2\omega(s-\alpha+1/4) \, \sin 2 \omega(s-\beta +1/4)}
{\sin 2\omega (s+1/4) \, \sin 2 \omega(s-1/4)},\\
C_s &=& \frac{\sin 2\omega(s+\alpha+1/4) \, \sin 2 \omega(s+\beta +1/4)}
{\sin 2\omega (s+1/4) \, \sin 2 \omega(s+3/4)},
\ea
$$
\lambda_n = - 4 \sin \omega n \, \sin \omega(n+2\alpha +2\beta -1),
$$
$$
x_s= 2 \cos 2\omega (s+1/4).
$$
Under the conditions \re{I} or \re{II}, this difference equation can be
transformed by a gauge transformation $P_n(x_s) = \chi(s) \psi_n(s)$ with
an appropriate periodic function $\chi(s)$ to the hermitian form
\be
U_{s+1} \psi_n(s+1) +U_s \psi_n(s-1) -(A_s +C_s) \psi_n(s) =\lambda_n
\psi_n(s),
\lab{herm2} \ee
where $U_s = \sqrt{A_s C_{s-1}}$ is a real function because $A_sC_{s-1}>0$
for all integer $s$. Hence \re{herm2} provides a nontrivial example of the
Hermitian Jacobi matrix with periodic coefficients whose
eigenfunctions and eigenvalues are known explicitly. Note that in this
case, the difference equation for {\it symmetric} polynomials is determined
by the {\it tridiagonal} Jacobi matrix (cf. \cite{Sz2}).

The weight function for the corresponding AWP has the form
\be
w_s = R \sin 2\omega(s+1/4) \prod_{k=1}^s {\frac{\sin 2 \omega(k+ \alpha
-3/4) \sin 2 \omega(k+ \beta - 3/4)}{\sin 2 \omega(k - \alpha +1/4)
\sin 2\omega(k - \beta +1/4)}},
\lab{wsymm} \ee
where $s=0, 1, \dots, N-1,$ and $R$ is a normalization constant.

One can check that $w_s>0$ for the intervals \re{I} and \re{II} as
it should be (by the Favard theorem) for the case of positive recurrence
coefficients $u_n$. As in the generic case, the considered symmetric
polynomials do not obey complete self-duality.

\section{Finite orthogonality relation for continuous $q$-Jacobi polynomials}
\setcounter{equation} 0

In this section we study a finite set of the continuous $q$-Jacobi
polynomials \cite{AW,GR,KS} for $q$ a root of unity. These polynomials appear
from the AWP with arbitrary parameters $a$ and $c$ and the choice
\be
b=aq^{1/2}, \quad d=cq^{1/2}
\lab{par1} \ee
(we use a slightly different parametrization than in
\cite{AW,GR,KS}). Moreover, we will assume that $M$ in \re{q}
{\it is odd}. Then from \re{EN} we
see that $E_N=0$, hence the zeros of the $N$th continuous
$q$-Jacobi polynomials are
\be x_s=2\cos
\pi {M\over N}(s+1/2)= q^{(s+1/2)/2}+q^{-(s+1/2)/2}, \quad s=0, 1, \dots, N-1.
\lab{rcj} \ee
We could set $M=1$ in \re{rcj}.
However this parametrization is more convenient for further calculations.
The zeros \re{rcj} are simple, real, and
they coincide with zeros of the monic Chebyshev polynomials of the first kind
$T_N(x)=2\cos{N\theta}$.

In order to simplify the expression for
$P_{N-1}(x_s)$, we use another form of the continuous Jacobi polynomials
due to Rahman \cite{AW,GR}
\be P_n(x_s)= D{_n} \; {{_4} \varphi_3}
\left({p^{-n},p^{-s},p^{s+1},-acp^{n}
\atop ap^{1/2},cp^{1/2}, -p} ;p, p\right),
\lab{aw2} \ee
where $p=q^{1/2}=e^{\pi i M/N}$ and
\be D_n=\frac{(-p; p)_n(ap^{1/2}; p)_n(cp^{1/2}; p)_n}
{p^{n/2}(-acp^n; p)_n}.
\lab{tilc} \ee
The formula \re{aw2} is valid for all $n=0,1,\dots,N-1$.
However, for $n=N-1$, a remarkable simplification takes place.  Since $M$ is
odd we have $p^{1-N}=-p$, hence \re{aw2} is reduced to the
balanced ${_3 \varphi_2}$ function
\be
 P_{N-1}(x_s)=   D_{N-1} \; {{_3} \varphi_2}
\left({p^{-s},p^{s+1},acp^{-1} \atop ap^{1/2},cp^{1/2}} ;p, p \right).
\lab{3f2} \ee
But this sum can be calculated exactly using Jackson's extension of the
Pfaff-Saalschutz formula \cite{GR} and adapting it to roots of unity.
After this, we have
\be P_{N-1}(x_s)= D_{N-1} \left( ac\over p\right)^s
\frac{(a^{-1}p^{3/2};p)_s(c^{-1}p^{3/2};p)_s}{(ap^{1/2};p)_s(cp^{1/2};p)_s}.
\lab{fff} \ee

The $P^{\prime}_N(x_s)$ can be easily calculated because
$P_N(x)=2\cos \theta N, \; x=2 \cos \theta$. Namely, we have
\be
P^{\prime}_N(x_s)= \frac{N(-1)^{s+(M-1)/2}}{\sin (\pi M(s+1/2)/N)}.
\lab{der} \ee
From \re{ws} we find thus the weight function of the
continuous $q$-Jacobi polynomials
\be
w_s =R_N  \sin (\pi M (s+1/2)/N)
\frac{(-p/ac)^s (ap^{1/2};p)_s (cp^{1/2};p)_s}{(a^{-1}p^{3/2};p)_s
(c^{-1}p^{3/2};p)_s}, \quad s=0, 1, \dots, N-1,
\lab{wcj} \ee
where the normalization constant is
\be
R_N=ip^{-1/2}
\frac{(-cp^{1/2};p)_{N-1}(-ap^{1/2};p)_{N-1}}{(-acp;p)_{N-1}(-p;p)_{N-1}}.
\lab{norm} \ee

In general, the parameters $a,c$ are complex and the weight function
\re{wcj} is also complex. Nevertheless, there are important cases when the
weight function becomes {\it positive}. Indeed, put $M=1$ and choose
the following parametrization \cite{AW}
\be a=p^{\alpha +1/2}, \quad c=-p^{\beta +1/2},
\lab{albet} \ee
which is quite natural for it leads to the ordinary Jacobi polynomials
$P_n^{(\alpha,\beta)}(x)$ when $N \to \infty$.
Using this parametrization and passing to the trigonometric notations,
we can state the following theorem.

{\bf Theorem 3.}
{\it The weight function of the finite-dimensional
$q$-analogs of the Jacobi polynomials \re{aw2} for $q=\exp(2\pi i/N)$
has the form
\be
w_s=R_N \sin(\omega(2s+1)) \prod_{k=1}^s
{\frac{\sin \omega(k+\alpha) \cos \omega(k+ \beta)}{\sin
\omega(k-\alpha) \cos \omega(k-\beta)}},
\lab{realw} \ee
where $\omega = \pi /2N$, $s=0, 1, \dots, N-1,$
and the normalization constant
\be
R_N=\prod_{k=1}^{N-1}{\frac{\sin\omega(k+\beta) \cos \omega(k+\alpha)}
{\cos\omega k \sin \omega(k+\alpha + \beta + 1 )}}.
\lab{realR} \ee
For
\be
-1 < \alpha,\, \beta < 1,
\lab{abc} \ee
the weight function \re{realw} is positive, $w_s>0.$
}
\medskip

Proof of the positivity of $w_s$ is straightforward.
In this case, the recurrence
coefficients $u_n$ are positive too and the coefficients $b_n$ are real. So the
condition \re{abc} leads to the class of ordinary orthogonal polynomials
\cite{Chi} having the positive weight function on the real spectral set. The
particular choice $\alpha=\beta$ corresponds to the $q$-ultraspherical
polynomials considered (for the case of $q$
a root of unity) in \cite{Sz2}. There are also nine special
cases when the weight function \re{realw} is reduced to simpler
expressions. These cases arise when $\alpha, \beta =-1/2, 0, 1/2$. For
example, when $\alpha=\beta=-1/2$ we get finite set of Chebyshev
polynomials of the first kind with the constant weight function $w_s=1/N$;
for $\alpha=\beta=0$ we obtain Legendre $q$-polynomials with $w_s=\sin\omega
\sin \omega(2s+1)$; for $\alpha=-1/2, \beta=1/2$ we have the weight function
$w_s=N^{-1}\cot\omega/2 \cos^2 \omega(s+1/2)$, etc.

Note that the normalization condition $w_s=1$ can be considered as a
nontrivial summation analog of the standard beta-integral, because in the
limit $q\to1$ (i.e., $N \to \infty$) we get the ordinary Jacobi polynomials
with continuous weight function $w(x)=(1-x)^{\alpha} (1+x)^{\beta}$.
Normalization of this function leads to the ordinary beta-integral.
$q$-Extension of the beta integral was found in \cite{AW}
when $q$ is inside the unit disk, $|q|<1$.

Consider some limiting cases. When $c \to 0$, we get the continuous
$q$-Laguerre polynomials \cite{KS} with the weight function
\be
w_s=ip^{-1/2} \frac{(-ap^{1/2};p)_{N-1}}{(-p;p)_{N-1}} \sin(\pi M
(s+1/2)/N)a^{-s} p^{-s^2/2} \frac{(ap^{1/2};p)_s}{(a^{-1}p^{3/2};p)_s}.
\lab{wcl} \ee
When $a=p^{1/2}$, we get an interesting trigonometric identity
$$
\sum_{s=0}^{N-1} \sin \pi \frac{M}{N}(s+1/2) \,p^{-s(s+1)/2}=-ip^{1/2}.
$$
If in \re{wcl} we take the limit $a \to 0$ then we get the weight
function for $q$-Hermite polynomials
\be
w_s=\frac{i(-1)^s p^{-s^2-s-1/2}}{(-p;p)_{N-1}} \sin( \pi M (s+1/2)/N).
\lab{wch} \ee

\section{Finite orthogonality for big $q$-Jacobi polynomials}
\setcounter{equation}0

Big $q$-Jacobi polynomials can be considered as some limiting subcase of the
AWP \cite{GR}. However, it is more convenient to work directly with these
polynomials themselves.
Their recurrence coefficients are
$$
u_n=\xi_n \eta_{n-1}, \quad b_n=1-\xi_n -\eta_n,
$$
where
\ba
\xi_n &=& -acq^{n+1}
\frac{(1-q^n)(1-bq^n)(1-abc^{-1}q^n)}{(1-abq^{2n})(1-abq^{2n+1})},
\nonumber \\
\eta_n &=& \frac{(1-aq^{n+1})(1-cq^{n+1})(1-abq^{n+1})}
{(1 - abq^{2n+1})(1-abq^{2n+2})}.
\lab{recbig} \ea
The explicit form of the big $q$-Jacobi polynomials is
\be
P_n(x)= D_n\; {{_3} \varphi_2} \left({q^{-n},abq^{n+1},x
\atop aq, cq} ;q, q \right),
\lab{big} \ee
where
$$D_n=\frac{(aq;q)_n(cq;q)_n}{(abq^{n+1};q)_n}.$$
For $q$ fixed by \re{q}, the $N$th polynomial has the form
\be
P_N(x)=x^N -1 + \frac{(1-a^N)(1-c^N)}{1-(ab)^N}.
\lab{PN3} \ee
In what follows, we set $c=1$. Then the zeros of the polynomial $P_N(x)$
coincide with the roots of unity and lie at the vertices of the regular polygon
\be
x_s=q^s, \quad s=1,2,\dots, N.
\lab{bigroot} \ee
For $P_{N-1}(x_s)$ we have
\be
P_{N-1}(x_s)=  D_{N-1} \; {{_2} \varphi_1}
\left({q^{s-N},ab \atop aq} ;q, q \right).
\lab{bigN-1} \ee
This expression can be evaluated using the root-of-unity analog
of the Chu-Vandermonde formula \cite{GR} yielding
\be
P_{N-1}(x_s)=D_{N-1} (ab)^{N-s}
\frac{(b^{-1}q; q)_{N-s}}{(aq; q)_{N-s}},
\lab{bigN-1,2} \ee
where
\be D_{N-1}=N\frac{(1-a^N)(1-abq^{-1})}{(1-a)(1-(ab)^N)}.
\lab{bigC} \ee
Hence for the weight function, using \re{ws}, we have
\be w_s= \frac{(1-a^N)(1-abq)(b;q)_s q^s}{aq(b-1)(1-a^Nb^N)(a^{-1};q)_s},
\quad s=1,2,\dots,N.
\lab{bigw} \ee
It can be verified that this weight function is normalized,
$\sum_{s=1}^N w_s=1.$

There are two interesting limiting cases: when $a \to \infty$ we get
the $q$-Meixner polynomials \cite{KS} with the weight function
$$w_s=b^{1-N}(1-b)^{-1}(b;q)_s q^s.$$
When $b \to 0$ we get the big
$q$-Laguerre polynomials \cite{KS} with the weight function
$$w_s=(qa)^{-1}(a^N-1) \frac{q^s}{(a^{-1};q)_s}.$$
It should be noted however
that in both these cases we obtain only special (one-parameter) type of the
corresponding polynomials due to the condition $c=1$.

\section{Finite orthogonality for alternative $q$-Jacobi polynomials}
\setcounter{equation}0

So far, we were considering AWP with the base $q$. In this section we consider
a special class of the AWP with the base $p=q^{1/2}=e^{\pi i M/N}$, where
$M$ is odd.
These polynomials are orthogonal on some points of the real interval and can
be interpreted as alternative $q$-Jacobi polynomials.

Let us set
$c=-d=1$ and keep $a$ and $b$ as arbitrary complex parameters for AWP.
Then the $n$th Askey-Wilson polynomial (with the
base $p$) is written as (cf. \re{aw})
\be
P_n(x)= D{_n}\:  {_4} \varphi_3 \left({p^{-n},-abp^{n-1},az,a/z \atop ab, a,
-a} ;p, p \right),
\lab{altjac}  \ee
where $x=z+z^{-1}$ and in the expression for normalization constant $D_n$
\re{C_n} we should substitute $p$ instead of $q$.
These polynomials do not belong to the family of finite polynomials
considered in Section 2, but rather they are subcases of the $q$-Racah
polynomials considered in \cite{AW0}. Still, as far as the authors understand,
the specific property of these polynomials to have
positive weight functions for some particular choices of $M$ and $N$
has not been discussed in the literature (see also \cite{Sz1,Sz2}).

Let $n=N+1$. Then taking into account that $p^N=-1$ we get
\be
P_{N+1}(x)= D_{N+1}\:  {_3} \varphi_2 \left({p^{-N-1},az,a/z \atop
a, -a} ;p, p \right),
\lab{ajacN+1} \ee
where
$$D_{N+1}=a^{-N-1}(1-a^2)(1-a^{2N}).$$
In order for polynomials $P_n(x), 0\leq n\leq N+1,$ to be well defined
we impose the constraints
\be
ab, a, b \neq p^k, \quad k=0, 1, \dots, 2N-1.
\lab{abrestr}\ee
In particular, in this case $u_n$ does not vanish between two zeros
at $n=0$ and $n=N+1$.

The expression \re{ajacN+1} can be simplified using the
Jackson summation formula \cite{GR}
\be
P_{N+1}(x)=(z;p)_{N+1}(z^{-1};p)_{N+1}= z^{-N-1}(z^2-1)(z^{2N}-1).
\lab{alN+1} \ee
From \re{alN+1} we find zeros $x_s$ of the polynomial $P_{N+1}(x)$
\be
x_s=p^s + p^{-s} =2\cos{\pi Ms/N}, \quad s=0,1,2,\dots, N.
\lab{alroots} \ee
It is seen from \re{alroots} that all the zeros $x_s$ are real and simple.

From \re{alN+1} one can derive
\be
P_{N+1}^{\prime}(x_s)=2N(-1)^s/\chi_s,
\lab{N+1'} \ee
where $\chi_s$ is a ``characteristic function" defined as follows
$$\chi_s=\cases{1, \qquad s=1,2,\dots,N-1 \cr 1/2, \quad s=0, N. \cr}$$

In order to calculate $P_N(x_s)$ note that AWP are invariant under any
permutation of the parameters $a,b,c,d$. Hence we can write
\be
P_{N}(x_s)= D_N \:  {_3} \varphi_2 \left({p^{-s},p^s,abp^{-1} \atop
a, b} ;p, p \right),
\lab{alN} \ee
where
$$D_N=\frac{(-1;p)_N(a;p)_N(b;p)_N}{(abp^{-1};p)_N}.$$
Again the expression \re{alN} can be simplified using the Jackson summation
formula:
\be
P_N(x_s)=D_N \frac{(p/a;p)_s (p/b;p)_s}{(a;p)_s(b;p)_s} (ab/p)^s, \quad
s=0,1,2,\dots, N.
\lab{alPNs} \ee
Combining all these results together we can state the following theorem.

{\bf Theorem 4.}
{\it The weight function of the finite set of alternative
$q$-Jacobi polynomials \re{altjac} with $p=\exp (\pi i M/N),$
$(N, M)=1, M$ odd,
and restrictions upon parameters \re{abrestr}, has the form
\be
w_s=\frac{h_N}{P_N(x_s)P_{N+1}^{\prime}(x_s)}=
A_N \, \chi_s \, \frac{(a;p)_s(b;p)_s(-p/ab)^s}{(p/a;p)_s(p/b;p)_s},
\lab{alw} \ee
where $s=0, 1, \dots, N$ and
$$A_N=\frac{(p;p)_N(-a;p)_N(-b;p)_N}{2N(-ab;p)_N}.$$
}
\medskip

The normalization condition $\sum_{s=0}^N{w_s}=1$ in this case leads to some
interesting trigonometric identities. Indeed, taking the limit $a \to 0, \:
b \to 0$ we obtain the weight function
\be
w_s=\chi_s \, \frac{(p;p)_N}{2N} \, (-1)^s \, p^{-s^2}.
\lab{a0b0} \ee
The normalization condition leads to the celebrated Gauss sum \cite{Chan}
(see also the Appendix)
\be \sum_{s=0}^{N-1}{(-1)^s\,
p^{-s^2}}=(-p;p)_{N-1}=\sqrt{N} p^{N(N-1)/4} \, \left({N\over M}\right),
\lab{gs1} \ee
where $({N \over M})$ is the Jacobi symbol. For $M=p_1p_2\cdots p_K$ --
the prime number decomposition of $M$ -- this symbol is defined
as the product
$$
\left({N \over M}\right)=\left({N \over p_1}\right)\cdots
\left({N \over p_K}\right),
$$
where $({N \over p_j})$ is the Legendre symbol. The latter is
equal to $+1$ if $N$ is a quadratic residue
(i.e., if there exists an integer $x$ such that $x^2=N\; (\mbox{mod}\; p_j)$),
to  $-1$ if $N$ is not a quadratic residue, and to $0$ if $p_j$ divides $N$.
In the simplest case $M=1$, one has $({N \over 1})=1$.

In another limiting case, setting $a=p^{1/2}$ and taking the limit $b \to 0$,
we get the weight function
\be
w_s=A_N\, \chi_s \, p^{-s^2/2} \lab{b0},
\ee
where
\be
A_N=\frac{(p;p)_N(-p^{1/2};p)_N}{2N}=\frac{(-p^{1/2};p)_{N}}{(-p;p)_{N-1}}.
\lab{Agauss} \ee
Using the formulas \re{pp} and \re{p1/2p} from the Appendix,
one can simplify the expression \re{Agauss} for $A_N$:
\be
A_N=\sqrt{\frac{2}{N}} \, p^{N/4} \left( {2N \over M} \right).
\lab{simplA} \ee
Hence we have the following trigonometric Gauss-like sum
\be
\sum_{s=0}^{N-1}{p^{-s^2/2}}={1\over 2}(1-
(-i)^{MN}) + \sqrt{\frac{N}{2}} \exp(-i \pi M/4) \left({2N \over M} \right).
\lab{newgauss} \ee
For $M=1$ this identity was derived in \cite{kor}. The $M\neq 1$
cases were probably discussed in the literature as well.

At the final stages of preparation of this paper, the authors have
discovered that in 1968 M. W. Wilson had considered a discrete analog of the
Legendre polynomials
\cite{wil}, which is a particular subcase of the above alternative
$q$-Jacobi polynomials. Indeed, let us set in \re{alw}
$a=-b=p^{1+\epsilon}$ and take the limit $\epsilon\to 0$. As a result,
one finds the weight function $w_0=w_N=0$ and
$w_s=\tan(\pi M/2N)\sin(\pi Ms/N)$ for $s=1, 2, \dots, N-1$.
For $M=1$, it coincides with the weight function of \cite{wil} up to a
normalization constant. Although some important characteristics of the
corresponding polynomials (e.g., the nice approximation properties)
have been found in \cite{wil} the description
was not complete. In particular, the explicit form of the recurrence
coefficients, which has not been derived, is
\begin{equation}
u_n={\sin^2\pi n/2N \cos \pi(n-1)/2N \cos \pi(n+1)/2N \over
\sin\pi(2n-1)/2N \sin\pi(2n+1)/2N}, \quad b_n=0.
\label{qlegwil}\end{equation}
Only the first three coefficients $u_{1, 2, 3}$ were found in \cite{wil}
and they coincide with \re{qlegwil}.
Note that the corresponding polynomials are $q$-analogs
of the Legendre polynomials different from those mentioned above
which were considered in \cite{Sz2}.
The latter are simpler but it is natural to expect that they provide
a similar approximation of the Legendre polynomials for $N\to\infty$.

Returning to the general weight function \re{alw}, we note that when
$a=p^{\alpha}, \; b=-p^{\beta}$, and $0 < \alpha, \beta <1$, the
recurrence coefficients $u_n$ become positive leading to the positivity of
the corresponding weight function $w_s, \; s=0,1,\dots,N$. Verification of this
statement is elementary.  In the limit $p \to 1$ we recover the ordinary
Jacobi polynomials $P_n^{(\alpha-1,\, \beta-1)}(x)$. We have thus
a subclass of the ``classical" orthogonal polynomials
(in the sense of positivity of the weight function on a real interval)
among the alternative $q$-Jacobi polynomials in the same way as for the
case of the continuous $q$-Jacobi polynomials.

Although we did not give a complete classification of the regions
of complex parameters $a, b, c, d$ and $q$ in the recurrence coefficients
\re{recaw} leading to positive measures for the corresponding (associated)
Askey-Wilson polynomials, the present work provides some necessary
tools for the resolution of this important problem.

The authors are deeply indebted to R. Askey and G. Gasper for helpful remarks
to the preliminary version of this paper.

\section{Appendix}
\setcounter{equation}0

In this appendix we present some useful formulas concerning $q$-calculus when
$q$ is a root of unity. Some of them are well known, nevertheless they are
included here because the systematic treatment of such formulas is
absent in the literature. We consider here only standard
$q$-hypergeometric series, although there are similar relations
for the bilateral series as well.

Assume that $q$ is a primitive root of unity, i.e., $q=e^{2\pi i M/N}$, where
$(N,M)=1$. Denote also $p=q^{1/2}=e^{\pi i M/N}$. The first set of
formulas is connected to the $q$-shifted factorial defined in the standard
manner
\be (a;q)_n=(1-a)(1-aq)\dots (1-aq^{n-1}),\quad (a;q)_0=1. \lab{qp} \ee
Obviously $(a;q)_{n+N}=(a;q)_n (a;q)_N$, so it is sufficient to study only the
case $0 \le n \le N$. The following identities are easy to derive
\be
(a;q)_N=1-a^N
\lab{aN} \ee
(in particular $(a;q)_N=0$ if $a$ is a primitive $N$th root of unity) and
\be
(a;q)_{N-1}=\frac{1-a^N}{1-aq^{-1}}, \quad \mbox{if $a\ne q$},
\lab{a,N-1} \ee
\be
(q;q)_{N-1}=N.
\lab{a,N} \ee
In general, there are no simple formulas for $(a;p)_N$ and $(a;p)_{N-1}$ like
\re{aN} and \re{a,N-1}. In our case the following relations for
inversion of the parameter are valid. Assume that $M$ is odd. Then
\ba
(a;p)_N=a^N p^{-N(N+1)/2} (-a^{-1}p;p)_N,
\lab{odd1} \\
(a;p)_{N-1}=a^{N-1} p^{1-N(N+1)/2} (-a^{-1}p^2;p)_{N-1}.
\lab{odd2} \ea
If $M$ is even (hence $N$ is odd), then
\ba
(a;p)_N = -a^N (a^{-1}p;p)_N,
\lab{even1} \\
(a;p)_{N-1} = pa^{N-1} (a^{-1}p^2;p)_{N-1}.
\lab{even2} \ea

The next set of formulas is connected to finite sums involving
the $q$-shifted factorial.
The simplest one is a special case of the $q$-binomial theorem which can be
written as
\be
{_1}\varphi_0(q^s;q, z)=\sum_{k=0}^{N-s}
{\frac{(q^{s};q)_k}{(q;q)_k}z^k}=\frac{1-z^N}{(z;q)_s}=(zq^{-1};q^{-1})_{N-s},
\lab{qbin} \ee
where $s=1,2,\dots,N-1$. In \re{qbin} and below we use the standard
notation for $q$-hypergeometric functions although we always assume
that only the first $N$ terms are present in the corresponding series.

The Chu-Vandermonde formula \cite{GR} for $q$ a root of unity looks like this:
\be
{{_2} \varphi_1}
\left({q^{s},a \atop c} ;q, q \right)=\frac{a^N-c^N}{1-c^N}\,
\frac{(q/c;q)_s}{(qa/c;q)_s}, \quad s=1,2,\dots,N-1,
\lab{Chu} \ee
where $a,c$ are arbirary complex parameters, $c^N \ne 1$. In particular, for
$s=1$ we get a pretty summation formula
\be
\sum_{k=0}^{N-1}{\frac{(a;q)_k}{(c;q)_k}q^k}=
\frac{(a^N-c^N)(1-cq^{-1})}{(1-c^N)(a-cq^{-1})}.
\lab{spChu} \ee
Similarly, Jackson's extension of the Pfaff-Saalschutz sum \cite{GR} in this
special case is written as
\be
\sum_{k=0}^{N-1}\frac{(a;q)_k (b;q)_k q^k}{(c;q)_k (abc^{-1}q^2;q)_k} =
\frac{(a^N-c^N)(b^N-c^N)(q-c)(abq-c)}{(1-c^N)(a^Nb^N-c^N)(qa-c)(qb-c)}.
\lab{saa} \ee
The Dixon $q$-summation formula \cite{GR} takes the form
\be
\sum_{k=0}^{N-1}{\frac{(-qa;q)_k(b;q)_k}{(-a;q)_k(a^2b^{-1}q;q)_k}(a/b)^k}=
\frac{(b-a^2)(1+a^N)(b^N-a^N)}{(b-a)(1+a)(b^N-a^{2N})}.
\lab{Dic} \ee
A special case of the Singh quadratic transformation yields, for $q^N=1$,
the summation formula
\be
{{_4} \varphi_3} \left({q,a,a/p,a^2q/b^2
\atop a^2, -ab^{-1}q, -ab^{-1}p^3} ;q, q \right)=
\frac{(-p;p)_{N-1}(b/p;p)_{N-1}}{(-a;p)_{N-1}(b/a;p)_{N-1}},
\lab{Sing} \ee
where we assume that $M$ is odd so that $p^N=-1$. If one takes the limit $a\to
0,\:  b\to 0, \: b/a \to 0$ then we get an important summation formula
\be
\sum_{k=0}^{N-1}(-1)^k q^{-k^2/2} = (-p;p)_{N-1},
\lab{gauss} \ee
going back to the works of Gauss in the beginning of nineteenth century.
The Gauss sum \re{gauss}
plays a prominent role in some problems of the number theory \cite{Chan}.
The right hand side of \re{gauss} can be rewritten as follows \cite{Chan}
\be
(-p;p)_{N-1} = \left({N \over M} \right) \sqrt{N} e^{\pi i M(N-1)/4},
\lab{pp} \ee
where $({N\over M})$ is the Jacobi symbol.
Analogously, considering the Gauss sum
\be
\sum_{k=0}^{2N-1}(-1)^k p^{-k^2/2} = (-p^{1/2};p^{1/2})_{2N-1},
\lab{gauss3} \ee
one can get the identity
\be
(-p^{1/2};p^{1/2})_{2N-1}=\left({2N \over M} \right) \sqrt{2N}
e^{\pi i M(2N-1)/4}. \lab{pp1/2} \ee
From this relation and the obvious equality
$$
(-p;p)_{N-1} \, (-p^{1/2};p)_N = (-p^{1/2};p^{1/2})_{2N-1}
$$
follows another interesting identity
\be
(-p^{1/2};p)_N=\frac{(-p^{1/2};p^{1/2})_{2N-1}}{(-p;p)_{N-1}} =
\left({2 \over M} \right) \: \sqrt{2} \, e^{\pi i MN/4} =
(-1)^{(M^2-1)/8}\: \sqrt{2} \: e^{i \pi MN/4},
\lab{p1/2p} \ee
where the following properties of the Jacobi symbol were used
$$
\left({m \over M} \right) \, \left({n \over M} \right) = \left({mn \over M}
\right), \quad \left({2 \over M} \right)= (-1)^{(M^2-1)/8}.
$$

\bb{99}

\bi{AAA} W. Al-Salam, W. R. Allaway, and R. Askey, {\it Sieved ultraspherical
polynomials}, Trans. Amer. Math. Soc. {\bf 284} (1984), 39-55.

\bi{AW0} R. Askey and J. Wilson, {\it A set of orthogonal polynomials that
generalize the Racah coefficients or $6-j$ symbols}, SIAM J. Math. Anal.
{\bf 10} (1979), 1008-1016.

\bi{AW} R. Askey and J. Wilson, {\it Some basic hypergeometric orthogonal
polynomials that generalize Jacobi polynomials}, Mem. Amer. Math. Soc.
{\bf 54} (1985), 1-55.

\bi{At} F. V. Atkinson, ``Discrete and continuous boundary problems''
in {\it Mathematics in Science and Engineering}, {\bf 8},
Academic Press, New York, 1964.

\bi{BI} E. Bannai and T. Ito, {\it Algebraic Combinatorics I. Association
Schemes}, Benjamin/Cummings, Menlo Park, 1984.

\bi{Car} R. D. Carmichael, {\it The general theory of linear $q$-difference
equations}, Amer. J. Math. {\bf 34} (1912), 147-168.

\bi{Chan} K. Chandrasekharan, {\it Introduction to Analytic Number Theory},
Grundlehren Math. Wiss. {\bf 148}, Springer-Verlag, New York, 1968.

\bi{Chi} T. S. Chihara, {\it An Introduction to Orthogonal Polynomials},
Math. Appl. {\bf 13}, Gordon and Breach, New York, 1978.

\bi{GR} G. Gasper and M. Rahman, {\it Basic Hypergeometric Series},
Encyclopedia Math. Appl. {\bf 35},
Cambridge University Press, Cambridge, 1990.

\bi{mas} D. Gupta and D. Masson, {\it Contiguous relations, continued
fractions and orthogonality}, to appear in Trans. Amer. Math. Soc.

\bi{IR} M. E. H. Ismail and M. Rahman, {\it The associated Askey-Wilson
polynomials}, Trans. Amer. Math. Soc. {\bf 328} (1991), 201-237.

\bi{KS} R. Koekoek and R. F. Swarttouw, {\it The Askey scheme of
hypergeometric orthogonal polynomials and its $q$-analog},
Delft University Technical Report 94-05, 1994.

\bi{kor} N. M. Korobov, {\it Exponential Sums and their Applications},
Math. Appl. (Soviet Ser.) {\bf 80}, Kluwer, Dordrecht, 1992.

\bi{Leo} D. A. Leonard, {\it Orthogonal polynomials, duality, and association
schemes}, SIAM J. Math. Anal. {\bf 13} (1982), 656-663.

\bi{NSU} A. F. Nikiforov, S. K. Suslov and V. B. Uvarov,
{\it Classical Orthogonal Polynomials of Discrete Variable}, Springer Ser.
Comput. Phys., Springer-Verlag, Berlin, 1991.

\bi{Sz1} V. Spiridonov and A. Zhedanov, {\it Discrete Darboux transformations,
the discrete time Toda lattice, and the Askey-Wilson polynomials},
Methods Appl. Anal. {\bf 2} (1995), 369-398.

\bi{Sz2} V. Spiridonov and A. Zhedanov, {\it $q$-Ultraspherical
polynomials for $q$ a root of unity}, Lett. Math. Phys. {\bf 37}
(1996), 173-180.

\bi{wil} M. W. Wilson,  {\it On a new discrete analogue of the Legendre
polynomials}, SIAM J. Math. Anal. {\bf 3} (1972), 157-169.

\eb

V. Spiridonov: Centre de Recherches Math\'ematiques, Universit\'e de Montr\'eal,
C.P. 6128, succ. Centre-ville, Montr\'eal (Qu\'ebec) H3C 3J7, Canada;
current address: Laboratory of Theoretical Physics, JINR, Dubna 141980
Moscow region, Russia
\medskip

A. Zhedanov: Donetsk Institute for Physics and Technology, Donetsk 340114, Ukraine

\end{document}